\documentstyle[11pt,twoside]{article}

\input{mssymb}
\makeatother

\newcommand{\mysec}[2]{%
\section*{\normalsize\hfil\sc {#1}. {#2}\hfill}%
\setcounter{theo}{0}\setcounter{equation}{0}\setcounter{section}{#1}%
\typeout{#1. #2}
\noindent}

\newcommand{\Proof}{\par\noindent{\em Proof. }}
\newcommand{\eop}{\nopagebreak\hspace*{\fill}$\Box$}

\newtheorem{theo}{Theorem}[section]

\newtheorem{cor}[theo]{Corollary}

\newtheorem{prop}[theo]{Proposition}
\newtheorem{definition}[theo]{Definition}

\newcounter{abc}   % Counter fr statements-environment wird deklariert
\newcounter{iiiii} % Counter fr aequivalenz-environment wird deklariert

\newenvironment{aequivalenz}
{\setcounter{iiiii}{0}
\begin{list}%
{{\rm (\roman{iiiii})}}%  Falls die items nicht angegeben sind: i)u.s.w.
{\usecounter{iiiii}
\parsep=0pt plus 1pt
\topsep=1pt plus 2pt minus 1pt
\itemsep=1pt plus 2pt minus 1pt
\leftmargin=3\baselineskip
\labelsep=.6\baselineskip
\labelwidth=2.4\baselineskip
\rightmargin 0pt}%
}%               Das war das zweite Argument von "newenvironment"
{\end{list}}

\newenvironment{statements}%
{\setcounter{abc}{0}
\begin{list}%
{{\rm (\alph{abc})}}%  Falls die items nicht angegeben sind: (a) u.s.w.
{\usecounter{abc}
\parsep=0pt plus 1pt
\topsep=1pt plus 2pt minus 1pt
\itemsep=1pt plus 2pt minus 1pt
\leftmargin=3\baselineskip
\labelsep=.6\baselineskip
\labelwidth=2.4\baselineskip
\rightmargin 0pt}%
}%               Das war das zweite Argument von "newenvironment"
{\end{list}}

{\begin{list}%
{$\bullet$}%  Falls die items nicht angegeben sind: Punkte
{\leftmargin=5\baselineskip
\labelsep=\baselineskip
\labelwidth=2.5\baselineskip
\parsep=0pt plus 1pt
\topsep=1pt plus 2pt minus 1pt
\itemsep=1pt plus 2pt minus 1pt
\rightmargin 0pt}%
}%               Das war das zweite Argument von "newenvironment"
{\end{list}}

\newif\ifrefsc
\let\thebibliographyalt=\thebibliography                                %
\def\thebibliography#1                                                  %
 {\ifrefsc
 \section*{\normalsize\hfil\sc References\hfill}
 \small
 \list                    %
 {[\arabic{enumi}]}{\settowidth\labelwidth{[#1]}\leftmargin\labelwidth  %
 \advance\leftmargin\labelsep                                           %
 \usecounter{enumi}}                                                    %
 \def\newblock{\hskip .11em plus .33em minus .07em}                     %
 \sloppy\clubpenalty4000\widowpenalty4000\vfuzz 2pt                      %
 \sfcode`\.=1000\relax                                                  %
 \else\thebibliographyalt{#1}\fi}                                         %

\makeatletter
\def\ps@smallheadings{            % Nur bei doppelseitigem Ausdruck
  \def\@oddfoot{}                 % Aufruf mit: pagestyle{smallheadings}
  \def\@evenfoot{}                % keine Footer
  \def\@evenhead{\hbox to \textwidth {%
  \vbox{\hbox to \textwidth 
        {\thepage\hfil\strut{\footnotesize\sc Kamil John and 
                                              Dirk Werner}\hfil}\vss}}}
  \def\@oddhead{\hbox to \textwidth {%
\vbox{\hbox to \textwidth 
        {\hfil{\footnotesize\sc\Kurztitel}\hfil\strut \thepage}\vss}}}
    }
\makeatother

\makeatletter
\def\eqalignno#1{\displ@y \tabskip\@centering
  \halign to\displaywidth{\hfil$\@lign\displaystyle{##}~$\tabskip\z@skip
    &$\@lign\displaystyle{{}##}$\hfil\tabskip\@centering
    &\llap{$\@lign##$}\tabskip\z@skip\crcr
    #1\crcr}}

\makeatother

\def\ersteSeite{\vspace*{32pt plus 2pt minus 2pt}\begin{center}
{\Large\sf\Titel}\\[15pt]{\sc\Autor}\\[26pt plus 2pt minus 2pt]
\end{center}\Abstrakt\vspace{0pt plus 2pt }\thispagestyle{empty}}

\def\Abstrakt{\begin{quote}\small\noindent{\sc Abstract.}
\Abstrakttext\end{quote}}

\newcommand{\N}{{\Bbb N}}
\newcommand{\T}{{\Bbb T}}
\newcommand{\Z}{{\Bbb Z}}

\newcommand{\Id}{I\mkern-1mud}

\newcommand{\eps}{\varepsilon}

\newcommand{\qqfa}{\qquad\forall}

\newcommand{\bea}{\begin{eqnarray*}}
\newcommand{\eea}{\end{eqnarray*}}
\newcommand{\beq}{\begin{equation}}
\newcommand{\eeq}{\end{equation}}
\newcommand{\begsta}{\begin{statements}}
\def\endsta{\end{statements}}
\newcommand{\begaeq}{\begin{aequivalenz}}
\def\endaeq{\end{aequivalenz}}

\newcommand{\dopu}{{:}\allowbreak\ }

\newcommand{\rest}[2]{#1\raisebox{-0.3ex}{\mbox{$\mid_{#2}$}}} 

\refsctrue

\begin{document}

%\preliminaryversion
%\vspace*{-10pt}

\def\Titel{ $M$-ideals of compact operators into $\ell_{p}$ }
\def\Kurztitel{\Titel     }
\def\Autor{Kamil John\footnote{supported by the grants 
of GA~AV\v CR No.~1019504 and of GA~\v CR No.~201/94/0069.} 
                      and Dirk Werner}
\def\Abstrakttext{
We show for $2\le p<\infty$ and subspaces $X$ of quotients of
$L_{p}$ with a $1$-unconditional 
finite-dimensional Schauder decomposition that $K(X,\ell_{p})$ is an
$M$-ideal in $L(X,\ell_{p})$.
}

\ersteSeite

\mysec{1}{Introduction}%
A closed subspace $J$ of a Banach space $X$ is called an $M$-ideal if
the dual space $X^{*}$ decomposes into an $\ell_{1}$-direct sum $X^{*}
= J^{\perp} \oplus_{1} V$, where 
$J^{\perp}  = \{x^*\in X^*\dopu \rest{x^*}{J}=0\}$ is the annihilator of
$J$ and $V$ is some closed subspace of~$X^{*}$. This notion is due to
Alfsen and Effros \cite{AlEf}, and it is studied in detail in
\cite{HWW}.

It has long been known that the space of compact operators
$K(\ell_{p})$ is an $M$-ideal in the space of bounded operators
$L(\ell_{p})$ for $1<p<\infty$ whereas this property fails for
$L_{p}=L_{p}[0,1]$ unless $p=2$; cf.\ Section~VI.4 in \cite{HWW}.
More recently, it was shown in \cite{KalW2} that $K(L_{p},\ell_{p})$
is an $M$-ideal if $1<  p \le2$, and it is not an $M$-ideal if
$p>2$. 

In this paper we wish to examine the $M$-ideal character of
$K(X,\ell_{p})$ for subspaces $X$  of
quotients of  $L_{p}$ and $2\le p<\infty$. Our idea
is to exploit the fact that those $X$ have Rademacher cotype~$p$ with
constant~1. This leads to the result mentioned in the abstract.

We would like to thank N.~Kalton and E.~Oja for their comments on
preliminary versions of this paper.

\mysec{2}{Results}%
Here is our main result.

\begin{theo} \label{3.7}
Let $1< p< \infty$ and suppose that the
Banach space $X$ admits a sequence of operators $K_n\in K(X)$
satisfying
\begsta
\item
$K_nx\to x$ for all $x\in X$,
\item 
$K^*_nx^*\to x^*$ for all $x^*\in X^*$,
\item
$\|\Id_X-2K_n\|\to 1$.
\endsta
Then $K(X,\ell_p)$ is
an $M$-ideal in $L(X,\ell_p)$ if 
\beq \label{eq7}
\limsup_n (\|x\|^{p}+\|x_n\|^{p})^{1/ p}\le
\limsup_n\biggl({\|x+x_n\|^{p}+\|x-x_n\|^{p}\over
2}\biggr)^{1/ p}
%              \eqno (7)
\eeq
for all $x,x_n\in X$ such that $x_n\to0$ weakly.
\end{theo}
\Proof 
Let $T\dopu X\to \ell_{p}$ be a contraction. We shall show that $T$
has property~$(M)$, i.e.,
$$
\limsup_{n} \|y+Tx_{n}\| \le \limsup_{n} \|x+x_{n}\|
$$
whenever $x\in X$, $y\in \ell_{p}$, $\|y\|\le\|x\|$, and $x_{n}\to 0$
weakly in $X$. This implies our claim by \cite[Th.~6.3]{KalW2}.

In fact, we have
\bea
\limsup_{n} \|y+Tx_{n}\|
&=&
\limsup_{n} \bigl( \|y\|^{p} + \|Tx_{n}\|^{p} \bigr)^{1/p} \\
&\le&
\limsup_{n} \bigl( \|x\|^{p} + \|x_{n}\|^{p} \bigr)^{1/p} \\
&\le&
\limsup_n\biggl({\|x+x_n\|^{p}+\|x-x_n\|^{p}\over
2}\biggr)^{1/ p} ;
\eea
so it is enough to show that
\beq\label{eq2.2}
\limsup_{n} \|x+x_{n}\|=
\limsup_{n} \|x-x_{n}\|.
\eeq
Let $\eps>0$. Pick $m\in \N$ so that
$$
\|K_{m}x-x\|\le\eps,\qquad \|\Id-2K_{m}\| \le 1+\eps.
$$
Then pick $n_{0}\in\N$ so that
$$
\|K_{m}x_{n}\| \le\eps \qqfa n\ge n_{0};
$$
this is possible since $x_{n}\to 0$ weakly and $K_{m}$ is compact. We
now have for $n\ge n_{0}$
\bea
(1+\eps) \|x_{n}+x\| &\ge&
\|(\Id-2K_{m})(x_{n}+x)\| \\
&=&
\| x_{n}-x - 2K_{m}x_{n} +2x -2K_{m}x\| \\
&\ge&
\|x_{n}-x\| -2\eps -2\eps
\eea
so that
$$
\limsup_{n} \|x_{n}+x\| \ge \limsup_{n} \|x_{n}-x\|,
$$
and by symmetry equality holds.
\eop

\bigskip
We note that (\ref{eq7}) is not a necessary condition, for
essentially trivial reasons: e.g., if $p<2$ and $X=\ell_{2}$, then
every operator from $X$ to $\ell_{p}$ is compact and, therefore,
$K(X,\ell_{p})$ is an $M$-ideal, but (\ref{eq7}) fails.

As the proof shows, one can as well consider all the Banach spaces sharing
the property
$$
\limsup_{n} \|y+y_{n}\| \le \limsup_{n}  \bigl( \|y\|^{p}+\|y_{n}\|^{p}
\bigr)^{1/p}
$$
whenever $y_{n} \to0 $ weakly, e.g., $\ell_{q}$ or the Lorentz spaces
$d(w,q)$ for $p\le q< \infty$. So our theorem is closely related to
\cite[Th.~3]{OjaCR} and \cite[Prop.~4.2]{Dirk5}.
Actually, we  needed assumptions (a)--(c) only to
ensure~(\ref{eq2.2}), a condition that could be called property~$(wM)$
in accordance with Lima's property~$(wM^{*})$ \cite{Lim95}.

Now we wish to give more concrete examples where Theorem~\ref{3.7}
applies. There is a natural class of Banach spaces in which
inequality~(\ref{eq7}) is valid. Recall that a Banach space $X$ has
Rademacher type~$p$ with constant~$C$ if for all
finite families $\{x_{1},\ldots, x_{n}\}\subset X$, with $r_{1},
r_{2}, \ldots$ denoting the Rademacher functions,
$$
\left( \int_{0}^{1}\,\biggl\| \sum_{k=1}^{n} r_{k}(t)x_{k} \biggr\|^{p}
       dt \right)^{1/p}
\le  C
\biggl( \sum_{k=1}^{n} \|x_{k}\|^{p} \biggr)^{1/p};
$$
it has Rademacher cotype~$p$ with constant~$C$ if
$$
\biggl( \sum_{k=1}^{n} \|x_{k}\|^{p} \biggr)^{1/p}
\le C
\left( \int_{0}^{1}\,\biggl\| \sum_{k=1}^{n} r_{k}(t)x_{k} \biggr\|^{p}
       dt \right)^{1/p}
$$
instead. Thus we see that 
the inequality (\ref{eq7}) is always satisfied when  $X$ has
Rademacher cotype~$p$  with  constant~1, which is the
case if $X$ is a subspace of a quotient of $L_{p}$ for $2\le p<\infty$. As
for assumptions (a)--(c) from Theorem~\ref{3.7}, these conditions are
obviously fulfilled if $X$ has a shrinking $1$-unconditional 
finite-dimensional Schauder decomposition or merely the shrinking
unconditional metric compact approximation property of \cite{CasKal}
and \cite{GKS}. 
Let us mention that the ``shrinking'' character of
these properties holds, by a well-known convex combinations argument
(cf.\ \cite[Lemma~VI.4.9]{HWW}), for reflexive spaces automatically.
These  observations yield the next corollary.

\begin{cor} \label{3.9}
Let $X$ be a subspace of a quotient of 
$L_p$, $2\le p<\infty$, and let $X$ have a  $1$-unconditional 
finite-dimensional Schauder decomposition or merely the 
unconditional metric compact approximation property.
Then $K(X,\ell_p)$ is an $M$-ideal in $L(X,\ell_p)$.
\end{cor}

More explicitly, we note that for instance $\ell_{p}$, $\ell_{p}
\oplus_{p} \ell_{r}$ and $\ell_{p}(\ell_{r})$, where $2\le r \le p
<\infty$,
satisfy these assumptions; but for these spaces the result of
Corollary~\ref{3.9} has already been known from \cite{Dirk5} or
\cite[p.~327]{HWW}. Yet there are other examples. In fact,
 Li \cite{Li-UCMAP} 
has exhibited spaces of
$\Lambda$-spectral functions $L^{p}_{\Lambda}(\T)$ for certain
$\Lambda\subset\Z$ that enjoy the 
unconditional metric compact approximation property.  Moreover, since
for $2\le q \le p<\infty$ the space $L_{q}$ is isometric to a quotient of
$L_{p}$, one can substitute ${q}$ for ${p}$ 
in the above list of examples.

Another way to see that (\ref{eq7})  holds for $L_{p}$, $2\le p
<\infty$, is
to observe that (\ref{eq7})  follows immediately from Clarkson's
inequality in $L_{p}$, that is
$$
\|f\|^{p} + \|g\|^{p} \le \frac{\|f+g\|^{p} + \|f-g\|^{p}}2
$$
for $p\ge 2$.
 Now, Clarkson's inequalities are valid in the
Schatten classes as well \cite{McCar}. Therefore we obtain a
noncommutative version of the previous corollary.
(Actually, this argument is not that different, because the Clarkson
inequality entails the desired cotype property.)

\begin{cor} \label{3.9a}
Let $X$ be a subspace of a quotient of the Schatten class
$c_p$, $2\le p<\infty$, and let $X$ have a  $1$-unconditional 
finite-dimensional Schauder decomposition or merely the 
unconditional metric compact approximation property.
Then $K(X,\ell_p)$ is an $M$-ideal in $L(X,\ell_p)$.
\end{cor}

There is a dual version of Theorem~\ref{3.7} which we state for
completeness.

\begin{theo} \label{3.3}
Let $1< p< \infty$ and $1/p + 1/p' = 1$. Suppose that the
Banach space $Y$ admits a sequence of operators $K_n\in K(Y)$
satisfying
\begsta
\item
$K_ny\to y$ for all $y\in Y$,
\item
$K^*_ny^*\to y^*$ for all $y^*\in Y^*$,
\item
$\|\Id_Y-2K_n\|\to 1$.
\endsta
Then $K(\ell_p,Y)$ is
an $M$-ideal in $L(\ell_p,Y)$ if 
\beq \label{eq6}
\limsup_n (\|y^*\|^{p'}+\|y_n^*\|^{p'})^{1/ p'}\le
\limsup_n\biggl({\|y^*+y_n^*\|^{p'}+\|y^*-y^*_n\|^{p'}\over
2}\biggr)^{1/ p'}
%            \eqno (6)
\eeq
for all $y^*,y^*_n\in Y^*$ such that $y^*_n\to0$ weak$^{\,*}$.
\end{theo}

The proof of Theorem~\ref{3.3} can be accomplished along the same
lines as above using property~$(M^{*})$ of a contraction 
(cf.\ \cite[p.~171]{KalW2}) instead.

Again, inequality (\ref{eq6}) is always satisfied when  $Y^*$ has
Rademacher cotype~$p'$  with  constant~1, which is the
case if $Y$ has Rademacher type~$p$ with constant~1. 
The latter holds
if $Y$ is a subspace of a quotient of $L_{p}$ or $c_{p}$ for $1<p\le2$.

\hfuzz 5pt

\mysec{3}{Concluding remarks}%
The conditions (\ref{eq7}) and (\ref{eq6}) can be understood as
averaging conditions. In an earlier draft of this manuscript we used
these conditions to establish what we call $p$-averaged versions of the
properties~$(M)$ and $(M^{*})$ of contractions~$T$, that is
$$
\limsup_n\|y + Tx_n\|
\le      
\cases { \displaystyle
\limsup_n\biggl({\|x+x_n\|^p+
\|x-x_n\|^p \over 2}\biggr)^{1/ p}
& for  $ p < \infty$ \cr
\displaystyle
\limsup_n\max(\|x+x_n\|,\|x-x_n\|)
& for $p=\infty$ \cr
}$$
whenever $x\in X$, $y\in Y$  with $\|y\|\le\|x\|$ and $x_{n}\to0$
weakly in $X$; respectively,
$$
\limsup_n\|x^*+T^*y^*_n\|
\le      
\cases { \displaystyle
\limsup_n\biggl({\|y^*+y_n^{*}\|^p+
\|y^*-y_n^{*}\|^p \over 2}\biggr)^{1/ p}
& for  $ p < \infty$ \cr
\displaystyle
\limsup_n\max(\|y^*+y_n^{*}\|,\|y^*-y_n^{*}\|)
& for $p=\infty$. \cr
}$$
for all $x^*\in X^*$, 
$y^*\in Y^*$ such that $ \|x^*\|\le \|y^*\|$ and
for all weak$^*$ null sequences
$(y^*_n)\subset Y^*$.
(As a matter of fact, (\ref{eq6}) implies the $p'$-averaged
property~$(M^{*})$ for a contraction
$T\dopu \ell_{p}\to Y$.)
 Using techniques from \cite{KalW2} (which in
turn depend on those from \cite{Kal-M})
one can prove the following results.

\begin{prop}\label{3.77}
Let $1\le p\le \infty$ and suppose that the
Banach space $X$ admits a sequence of operators $K_n\in K(X)$
satisfying
\begsta
\item
$K_nx\to x$ for all $x\in X$,
\item 
$K^*_nx^*\to x^*$ for all $x^*\in X^*$,
\item
$\|\Id_X-2K_n\|\to 1$.
\endsta
Let $Y$ be a Banach space. Then $K(X,Y)$ is an $M$-ideal in
$L(X,Y)$ if and only if every contraction $T\dopu X\to Y$ has 
$p$-averaged~$(M)$.
\end{prop}

\begin{prop}\label{3.73}
Let $1\le p\le \infty$ and suppose that the
Banach space $Y$ admits a sequence of operators $K_n\in K(Y)$
satisfying
\begsta
\item
$K_ny\to y$ for all $y\in Y$,
\item
$K^*_ny^*\to y^*$ for all $y^*\in Y^*$,
\item
$\|\Id_Y-2K_n\|\to 1$.
\endsta
Let $X$ be a Banach space. Then $K(X,Y)$ is an $M$-ideal in
$L(X,Y)$ if and only if every contraction $T \dopu X\to Y$ has
$p$-averaged~$(M^{*})$.
\end{prop}

It is well known (cf.\ \cite[Th.~I.2.2]{HWW})
that a closed subspace $J$ of a Banach space $X$ is an
$M$-ideal in $X$ if and only if the following
3-ball property holds: For all
$y_1,y_2,y_3 \in B_J$, all $x\in B_X$  and all
$\eps>0 $ there is $y\in J$ such that $\|x+y_i-y\|\le
1+\eps$ for $i=1,2,3$. (Here $B_{X}$ denotes the closed unit ball
of~$X$.) Upon replacing the number~3 by some $n\in \N$ we obtain the
$n$-ball property, which is equivalent to the 3-ball property provided
$n\ge3$.
One may ``average'' this condition as well and obtain the following
characterisation of $M$-ideals by means of an averaged 3-ball property.

\begin{prop}\label{2.1}
A closed subspace $J$ of a Banach space
$X$ is an $M$-ideal in $X$ if and only if 
\begsta
\item[\rm (A)]         \it
For all
$y_1,y_2,y_3 \in B_J$, $x\in B_X$ and 
$\eps>0 $ there is $y\in J$ such that 
$$
\|x+y_i-y\| +\|x-y_i-y\|\le 2(1+\eps) \quad {\it for\ }i=1,2,3.
$$
\endsta
holds.
\end{prop}
\Proof
Evidently the 6-ball property implies (A).
Conversely, suppose (A). In order to show that $J$ is an $M$-ideal in $X$
we will verify the ordinary 3-ball property (see above).
Now an inspection of the proof of
\cite[Theorem~I.2.2]{HWW} shows that one may additionally assume that
${\rm dist}(x,J)\ge 1-\eps$, in which case (A) implies that
$$
\|x+y_{i}-y\| \le 2(1+\eps) - \|x-y_{i}-y\| \le 1+3\eps,\qquad 
i=1,2,3,
$$
and we are done.
\eop

%
%
%  References
%
%
\typeout{References}
%\bibliography{az}

\begin{thebibliography}{1}

\bibitem{AlEf}
{\sc E.~M. Alfsen and E.~G. Effros}.
\newblock {\em Structure in real {Banach} spaces. {Parts} {I and II}}.
\newblock Ann. of Math.
\newblock {\bf 96} (1972), 98--173.

\bibitem{CasKal}
{\sc P.~G. Casazza and N.~J. Kalton}.
\newblock {\em Notes on approximation properties in separable {Banach} spaces}.
\newblock In: P.~F.~X. M{\"{u}}ller and W.~Schachermayer, editors, {\em
  Geometry of Banach Spaces, Proc. Conf. Strobl 1989, {\em London Mathematical
  Society Lecture Note Series 158}}, pages 49--63. Cambridge University Press,
  1990.

\bibitem{GKS}
{\sc G.~Godefroy, N.~J. Kalton, and P.~D. Saphar}.
\newblock {\em Unconditional ideals in {Banach} spaces}.
\newblock Studia Math.
\newblock {\bf 104} (1993), 13--59.

\bibitem{HWW}
{\sc P.~Harmand, D.~Werner, and W.~Werner}.
\newblock {\em {$M$}-{Ideals} in {Banach} {Spaces} and {Banach} {Algebras}}.
\newblock Lecture Notes in Math.\ 1547. Springer, Berlin-Heidelberg-New York,
  1993.

\bibitem{Kal-M}
{\sc N.~J. Kalton}.
\newblock {\em {$M$}-ideals of compact operators}.
\newblock Illinois J. Math.
\newblock {\bf 37} (1993), 147--169.

\bibitem{KalW2}
{\sc N.~J. Kalton and D.~Werner}.
\newblock {\em Property {$(M)$}, {$M$}-ideals and almost isometric structure of
  {Banach} spaces}.
\newblock J. Reine Angew. Math.
\newblock {\bf 461} (1995), 137--178.

\bibitem{Lim95}
{\sc \AA.~Lima}.
\newblock {\em Property {$(wM^*)$} and the unconditional 
                     metric compact approximation property}.
\newblock Studia Math.
\newblock {\bf 113} (1995), 249--263.

\bibitem{Li-UCMAP}
{\sc D.~Li}.
\newblock {\em Complex unconditional metric approximation property for
  {$C_\Lambda(\T)$} spaces}.
\newblock Preprint, 1995.

\bibitem{McCar}
{\sc Ch.~A. McCarthy}.
\newblock {\em $c_p$}.
\newblock Israel J. Math.
\newblock {\bf 5} (1967), 249--271.

\bibitem{OjaCR}
{\sc E.~Oja}.
\newblock {\em Dual de l'espace des op\'{e}rateurs lin\'{e}aires continus}.
\newblock C. R. Acad. Sc. Paris, S\'er. A
\newblock {\bf 309} (1989), 983--986.

\bibitem{Dirk5}
{\sc D.~Werner}.
\newblock {\em New classes of {B}anach spaces which are {$M$}-ideals in their
  bi\-duals}.
\newblock Math. Proc. Cambridge Phil. Soc.
\newblock {\bf 111} (1992), 337--354.

\end{thebibliography}
%\bibliographystyle{standard}

%             %
%
% Address
%
%
\small
\bigskip
\noindent
Mathematical Institute, Czech Academy of Sciences, \v{Z}itn\'{a} 25,
\\
CZ-11\,567 Prague 1, Czech Republic; \ 
e-mail:
kjohn@earn.cvut.cz

\smallskip\noindent
I.~Mathematisches Institut, Freie Universit\"at Berlin,
Arnimallee 2--6, \\
D-14\,195 Berlin, Germany; \ 
e-mail: 
werner@math.fu-berlin.de

\end{document}